\RequirePackage{ifpdf}
\ifpdf 
\documentclass[pdftex]{sigma}
\else
\documentclass{sigma}
\fi

\newcommand{\mg}{\mathfrak{g}}
\newcommand{\bul}{{\bullet}}

\renewcommand{\S}{{\rm S}}
\newcommand{\U}{{\rm U}}

\begin{document}

\allowdisplaybreaks
\renewcommand{\textfraction}{0.01}
\renewcommand{\topfraction}{0.99}

\renewcommand{\thefootnote}{$\star$}

\renewcommand{\PaperNumber}{060}

\FirstPageHeading

\ShortArticleName{Shoikhet's Conjecture and Duf\/lo Isomorphism on (Co)Invariants}

\ArticleName{Shoikhet's Conjecture and Duf\/lo Isomorphism \\ on (Co)Invariants\footnote{This paper is a
contribution to the Special Issue on Deformation Quantization. The
full collection is available at
\href{http://www.emis.de/journals/SIGMA/Deformation_Quantization.html}{http://www.emis.de/journals/SIGMA/Deformation\_{}Quantization.html}}}

\Author{Damien CALAQUE and Carlo A. ROSSI}

\AuthorNameForHeading{D. Calaque and C.A. Rossi}

\Address{Department of Mathematics, ETH Zurich, 8092 Zurich, Switzerland}

\Email{\href{mailto:damien.calaque@math.ethz.ch}{damien.calaque@math.ethz.ch}, \href{mailto:carlo.rossi@math.ethz.ch}{carlo.rossi@math.ethz.ch}}

\ArticleDates{Received May 23, 2008, in f\/inal form August 29,
2008; Published online September 03, 2008}

\Abstract{In this paper we prove a conjecture of B.~Shoikhet. This conjecture states that the tangent isomorphism on homology,
between the Poisson homology associated to a~Poisson structure on $\mathbb{R}^d$ and the Hochschild homology of its
quantized star-product algebra, is an isomorphism of modules over the (isomorphic) respective cohomology algebras.
As a~consequence, we obtain a version of the Duf\/lo isomorphism on coinvariants.}


\Classification{16E45; 16E40; 81Q30}

\renewcommand{\thefootnote}{\arabic{footnote}}
\setcounter{footnote}{0}

\section{Introduction}

In his seminal paper \cite{K} on the deformation quantization of Poisson manifolds M.~Kontsevich proved that
the dif\/ferential graded Lie algebra (shortly, DGLA) of polydif\/ferential operators on a smooth manifold $M$ is formal
(i.e.~it is quasi-isomorphic, as a DGLA, to its cohomology). As a~consequence one obtains that any Poisson structure
$\pi$ on $M$ can be quantized in the sense of~\cite{BFFLS}, and that the Hochschild cohomology of the deformed quantized
algebra is isomorphic to the Poisson cohomology of $(M,\pi)$. Moreover, it is known that the Hochschild cohomology of
an associative algebra is naturally equipped with an associative cup-product; and Kontsevich proved that the mentioned
isomorphism between Hochschild and Poisson cohomologies is actually multiplicative if $M=\mathbb{R}^d$, using a homotopy
argument involving the so-called {\em Kontsevich eye}.
The proof of this statement, known as the ``compatibility with cup-products'', has been clari\-f\/ied in~\cite{MT},
and appeared to have a surprising application to Lie theory~\cite{K} (see also~\cite{PT} and~\cite{BCKT}) in
providing a new proof, together with a cohomological extension, of the famous Duf\/lo isomorphism~\cite{Du}.
We recall to the reader that the result of Duf\/lo states that the Poincar\'e--Birkhof\/f--Witt map can be modif\/ied
so that it reduces to an isomorphism of algebras on invariants.

A homological version of Kontsevich's formality theorem has been formulated by B.~Tsygan in~\cite{Ts} and proved by Shoikhet in \cite{Sh}
(in the case $M=\mathbb{R}^d$) and by Tsygan and Tamarkin in~\cite{TT}, using dif\/ferent approaches.
It broadly states that the Hochschild chain complex of
the algebra of smooth functions on $M$ is formal as a DG Lie module over the DGLA of polydif\/ferential operators.
Again, one obtains as a direct consequence of the general formalism on $L_\infty$-algebras, their modules and
$L_\infty$-morphisms between them, that the Hochschild homology of the deformed quantized algebra is isomorphic
to the Poisson homology of $(M,\pi)$.
The present paper is mainly concerned about the multiplicativity of this isomorphism. Namely, Hochschild (resp.~Poisson)
homology is naturally a graded module over the Hochschild (resp.~Poisson) cohomology algeb\-ra, and we prove that
the isomorphism induced by the Tsygan--Shoikhet formality intertwines these module structures (both cohomology
algebras being themselves isomorphic thanks to the compatibility with cup-products).
We call our result the ``compatibility with cap-products''.

As in the cohomological situation, this compatibility with cap-products has an application to Lie theory in providing
a version of the Duf\/lo isomorphism for coinvariants.

The main goal of this paper is to present a short and comprehensible proof of this conjecture (which is a particular
case of a more general result, whose detailed proof is presented in \cite{CR1}). Here the proof also relies on a
homotopy argument, Kontsevich's eye being replaced by the {\em I-cube}, this time being a manifold with corners
of dimension 3.

The paper is organized as follows.
In Section~\ref{s-1} we state the main result we mentioned in this introduction.
Section~\ref{s-2} is a brief reminder on Kontsevich's and Shoikhet's conf\/iguration spaces and their compactif\/ications; we also mention how they are related to each other, a fact which will play a central r\^ole in some later computations.
In Section~\ref{s-3} we recall the construction of Kontsevich's and Shoikhet's formality $L_\infty$-quasi-isomorphisms. This is the f\/irst time where the conf\/iguration spaces that were introduced in the previous section appear operatively.
Section~\ref{s-4} is the heart of the paper.
It contains the proof of the main result, for which compactif\/ied conf\/iguration spaces (and integrals
over them) play a crucial r\^ole.
We f\/inally end the paper with the proof of a version of the Duf\/lo isomorphism for coinvariants, which we obtain as a consequence of our main result; we only observe that we prefer to give a direct computational proof, as opposed to the proof of the same result in~\cite{Sh}, where it was proved under the assumption that the conjecture (whose proof is the core of this paper) were true.

\section{The main result}\label{s-1}

For the manifold $V=\mathbb R^d$, we consider the dif\/ferential graded Lie algebras (shortly, DGLA)
$T_{\mathrm{poly}}^\bul(V)$ and $\mathcal D_{\mathrm{poly}}^\bul(V)$ of polyvector vector f\/ields on
$V$ and of polydif\/ferential operators on $V$ respectively.

Further, let $\gamma$ be a solution of the Maurer--Cartan equation in $T_{\mathrm{poly}}^\bul(V)$ of the form
\begin{equation}\label{eq-MCel}
\gamma=\hbar \pi,
\end{equation}
where $\pi$ is a bivector f\/ield and $\hbar$ a formal parameter (therefore in particular $\pi$ is a Poisson structure).
The Formality Theorem of Kontsevich~\cite{K} implies that the polydif\/ferential operator
\[
B=\mathcal U(\gamma)=\sum_{n\geq 1} \frac{1}{n!}\mathcal U_n(\underset{n}{\underbrace{\gamma,\dots,\gamma}})
\]
$\mathcal U_n$ being the Taylor components of the $L_\infty$-quasi isomorphism between $T_{\mathrm{poly}}^\bul(V)$ and
$\mathcal D_{\mathrm{poly}}^\bul(V)$, satisf\/ies the Maurer--Cartan equation in the DGLA
$\mathcal D_\mathrm{poly}^\bul(V)[\![\hbar]\!]$ of (series of) polydif\/ferential operators on $V$, viewed as a
subcomplex of the Hochschild cochain complex of the $\mathbb C[\![\hbar]\!]$-algebra  $A=C^\infty(V)[\![\hbar]\!]$.

We denote by $\mu$ the usual multiplication on the algebra $A$: it may be viewed as an element of
$\mathcal D_{\mathrm{poly}}^\bul(V)[\![\hbar]\!]$ of degree $1$, and the sum $\mu+B$ specif\/ies an associative
product $\star$ on $A$, which is a deformation of the usual product on $A$, and which moreover satisf\/ies
\[
f\star g-g\star f=2\hbar\langle\pi,\mathrm d f\wedge\mathrm d g\rangle+\mathcal O(\hbar^2).
\]
In other words, $\star$ is a quantization of the Poisson structure $\pi$ on $V$ in the sense of~\cite{BFFLS}.

The DGLA $T_{\mathrm{poly}}^\bul(V)$ possesses an associative product $\cup$, namely the usual $\wedge$-product on
$T_{\mathrm{poly}}^\bullet(V)$, and a solution $\gamma$ of the Maurer--Cartan equation def\/ines, by means of the
Schouten--Nijenhuis bracket\footnote{We observe that, if we follow the sign conventions of~\cite{AMM}, we have
then to modify the Schouten--Nijenhuis bracket as
$[\alpha,\beta]'_\mathrm{SN}=-[\beta,\alpha]_\mathrm{SN}$,
where the Schouten--Nijenhuis bracket on the right-hand side is the usual one.}, a dif\/ferential $\gamma\cdot=[\gamma,\ ]$ on
$T_{\mathrm{poly}}^\bul(V)$ w.r.t.\ $\cup$.

The (negatively graded) deRham complex $\Omega^{-\bul}(V)$ of dif\/ferential forms on $V$ is naturally a~dif\/ferential graded module (shortly, DGM)
over the DGLA $T_{\mathrm{poly}}^\bul(V)$: the extensions of the Lie derivative $\mathrm L$ by means of the Cartan formula and of the contraction operator $\iota$ def\/ine
respectively a dif\/ferential $\mathrm L_\gamma$, for $\gamma$ as in (\ref{eq-MCel}), and an action $\cap$ of $T_{\mathrm{poly}}^\bul(V)$ on
$\Omega^{-\bul}(V)$.

On the other hand, for $\gamma$ as above, there is a product $\cup$ (of degree $1$) on the (shifted by $1$) Hochschild cochain complex
$C^\bul(A,A)[1]$ of $A$ with values in $A$,
\begin{equation*}
\left(\varphi\cup\psi\right)\!(a_1,\dots,a_{p+q})=\varphi(a_1,\dots,a_p)\star\psi(a_{p+1},\dots,a_{p+q}).
\end{equation*}
Additionally, the Hochschild dif\/ferential $\mathrm d_\mathrm H$ on the Hochschild cochain complex of $A$ is modif\/ied to the Hochschild dif\/ferential
$\mathrm d_{\mathrm H,\star}$ w.r.t.\ $\star$.
All these structures descend to the subcomplex of polydif\/ferential operators on $V$.

For the algebra $A$, we consider the (negatively graded) Hochschild chain complex $C_{-\bullet}(A,A)$.
For $\gamma$ as (\ref{eq-MCel}), there is an action of $C^\bul(A,A)[1]$ on $C_{-\bullet}(A,A)$ via
\begin{equation*}
\varphi\cap (a_0|a_1|\cdots|a_n)=(a_0\star\varphi(a_1,\dots,a_m)|a_{m+1}|\cdots|a_n),
\end{equation*}
if $m\leq n$; if $m>n$, the action is trivial.
Furthermore, we also have the dif\/ferential $\mathrm b_\star$ on $C_{-\bullet}(A,A)$, which modif\/ies the usual Hochschild dif\/ferential $\mathrm b$ on
$C_{-\bul}(A,A)$.
The previous formula def\/ines also an action of the DGLA of polydif\/ferential operators on $V$ on $C_{-\bul}(A,A)$.

We denote by $\mathcal S_n$, $n\geq 0$, the Taylor components of the $L_\infty$-quasi-isomorphism $\mathcal S$ from the $L_\infty$-module
$C_{-\bullet}(A,A)$ to the DGM $\Omega^\bul (V)[\![\hbar]\!]$, both over $T_{\mathrm{poly}}^\bul(V)$, constructed in~\cite{Sh}.
The DGM structure of $\Omega^\bul (V)[\![\hbar]\!]$ over $T_{\mathrm{poly}}^\bul(V)$ comes from the Lie derivative of polyvector f\/ields
on dif\/ferential forms, while, as shown in~\cite{Sh}, composition of the $L_\infty$-quasi-isomorphism $\mathcal U$ with the action $\mathrm L$
of Hochschild cochains on $A$ on Hochschild chains gives $C_{-\bullet}(A,A)$ the structure of an $L_\infty$-module over $T_{\mathrm{poly}}^\bullet(V)$.

For a solution $\gamma$ of the Maurer--Cartan equation as in (\ref{eq-MCel}), we consider the following linear maps:
\begin{align}
\label{eq-tangcoh}\mathcal U_\gamma(\alpha)&=\sum_{n\geq 0}\frac{1}{n!} \mathcal U_{n+1}(\alpha,\underset{n}{\underbrace{\gamma,\dots,\gamma}}),\quad\text{resp.}\\
\label{eq-tangh}\mathcal S_\gamma(c)&=\sum_{n\geq 0}\frac{1}{n!} \mathcal S_n(\underset{n}{\underbrace{\gamma,\dots,\gamma}},c),
\end{align}
for a general polyvector f\/ield $\alpha$ on $V$, resp.\ Hochschild chain $c$ on $A$.

\begin{theorem}\label{t-cupprod}
For a solution $\gamma$ of the Maurer--Cartan equation as in \eqref{eq-MCel}, \eqref{eq-tangcoh} is a
quasi-isomor\-phism of complexes
\[
\mathcal U_\gamma:\left(T_{\mathrm{poly}}^\bul (V)[\![\hbar]\!],\gamma\cdot\right)\to \left(\mathcal D_{\mathrm{poly}}^\bul(V)[\![\hbar]\!],\mathrm d_{\mathrm H,\star}\right),
\]
which additionally preserves the products in the corresponding cohomologies,
\begin{equation*}
\left[\mathcal U_{\gamma}([\alpha]\cup[\beta])\right]=\left[\mathcal U_{\gamma}([\alpha])\cup\mathcal U_{\gamma}([\beta])\right],
\end{equation*}
square brackets denoting cohomology classes.
\end{theorem}
We refer to~\cite{K,MT} for the proof of Theorem~\ref{t-cupprod}.
The main result of this paper is the proof of Conjecture 3.5.3.1 in~\cite{Sh}, which we may state in the following
\begin{theorem}\label{t-capprod}
For a solution $\gamma$ of the Maurer--Cartan equation as in \eqref{eq-MCel}, \eqref{eq-tangh} is a quasi-isomor\-phism of complexes
\[
\mathcal S_\gamma:\left(C_{-\bullet}(A,A),\mathrm b_\star\right)\to\left(\Omega^\bul(V)[\![\hbar]\!],\mathrm L_\gamma \right)
\]
and additionally preserves the action of $T_{\mathrm{poly}}^\bul(V)$ in the corresponding cohomologies,
\begin{equation*}
\left[[\alpha]\cap \mathcal S_{\gamma}([c])\right]=\left[\mathcal S_{\gamma}(\mathcal U_{\gamma}([\alpha])\cap[c])\right],
\end{equation*}
with the previous notation for cohomology classes.
\end{theorem}

\section[Configuration spaces and their compactifications]{Conf\/iguration spaces and their compactif\/ications}\label{s-2}

We brief\/ly discuss in this Section conf\/iguration spaces of $i)$ points in the complex upper-half plane $\mathcal H$ and on the real axis
$\mathbb R$, and $ii)$ points in the interior of the punctured unit disk $D$ and on the boundary $S^1$, and their compactif\/ications \`a la
Fulton--MacPherson.

\subsection[Configuration spaces of points in the upper half-plane]{Conf\/iguration spaces of points in the upper half-plane}\label{ss-2-1}

For a pair of non-negative integers $(n,m)$, the (open) conf\/iguration space $C_{n,m}^+$ is def\/ined as
\begin{equation*}
C_{n,m}^+=\left\{(p_1,\dots,p_n,q_1,\dots,q_m)\in \mathcal H^n\times \mathbb R^m
:\ p_i\neq p_j,\ i\neq j,\ q_1<\cdots<q_m\right\}/G_2,
\end{equation*}
where $G_2$ is the semidirect product $\mathbb R^+\ltimes \mathbb R$, acting via rescalings and translations.
If  $2n + m-2 \geq  0$, $C_{n,m}^+$ is a smooth real manifold of dimension $2n+m-2$.
We may consider more general conf\/iguration spaces $C_{A,B}^+$, where $A$ is any f\/inite set and $B$ is any ordered
f\/inite set.

The conf\/iguration space $C_n$ is def\/ined as
\begin{equation*}
C_n=\left\{(p_1,\dots,p_n)\in \mathbb C^n:\ p_i\neq p_j,\ i\neq j\right\}/G_3,
\end{equation*}
where $G_3$ is the semidirect product $\mathbb R^+\ltimes \mathbb C$, acting via rescalings and complex translations.
If $2n-3\geq 0$, $C_n$ is a smooth real manifold of dimension $2n-3$.
Again, we may consider more general conf\/iguration spaces $C_A$, for any f\/inite set $A$.

Both conf\/iguration spaces $C_{A,B}^+$ and $C_A$ are orientable, see e.g.~\cite{AMM}.

Conf\/iguration spaces $C_{A,B}^+$ and $C_A$ admit compactif\/ications \`a la Fulton--MacPherson, denoted by
$\mathcal C_{A,B}^+$ and $\mathcal C_A$ respectively: they are smooth manifolds with corners and we refer
to~\cite{K, MT, CR} and~\cite{CR1} for their explicit constructions.

\subsection[Configuration spaces of points in the punctured disk]{Conf\/iguration spaces of points in the punctured disk}\label{ss-2-2}

As in Subsection~\ref{ss-2-1}, for a pair of non-negative integers $(n,m)$, $m\geq 1$, the (open) conf\/iguration
space $D_{n,m}^+$ is def\/ined as
\begin{gather*}
D_{n,m}^+=\left\{(p_1,\dots,p_n,q_1,\dots,q_m) \in (D^\times)^n \times  (S^1)^m:\!\begin{cases}
p_i\neq p_j, \!\! & \!\!\!\!  i\neq j,\\
q_1<\cdots<q_m<q_1,&
\end{cases}
\right\}\!/S^1,
\end{gather*}
where $D^\times$ denotes the punctured unit disk, and where we introduced a cyclic order on $S^1$; the group $S^1$
acts by rotations.
If $2n+m-1\geq 0$, $D_{n,m}^+$ is a smooth real manifold of dimension $2n+m-1$.
As before, we may consider conf\/iguration spaces $D_{A,B}^+$, where $A$ is any f\/inite set and $B$ is any
cyclically ordered f\/inite set.
When $|B|=1$, we omit the superscript $+$.

For a positive integer $n$, we consider the conf\/iguration space
\begin{equation*}
D_n=\left\{(p_1,\dots,p_n)\in (\mathbb C^\times)^n:\ p_i\neq p_j,\ i\neq j\right\}/\mathbb R^+,
\end{equation*}
where $\mathbb R^+$ acts by rescaling.
It is obviously a smooth real manifold of dimension $2n-1$, when $2n-1\geq 0$.
We may consider conf\/iguration spaces $D_A$, with $A$ any f\/inite set.

Finally, $D_{A,B}^+$ and $D_A$ are orientable, by the same arguments as in~\cite{AMM}.

Conf\/iguration spaces $D_{A,B}^+$ and $D_A$, admit compactif\/ications \`a la Fulton--MacPherson, denoted by
$\mathcal D_{A,B}^+$ and $\mathcal D_A$ respectively, which are smooth manifolds with corners.

Being $\mathcal D_{A,B}^+$ a stratif\/ied space, its boundary strata of codimension $1$ are given in the following list:
\begin{enumerate}\itemsep=0pt
\item[$i)$] There is a subset $A_1$ of $A$, obeying $1\leq |A_1|\leq |A|$, such that
\begin{equation}\label{eq-diskb1}
\partial_{A_1,0}\mathcal D_{A,B}^+\cong \mathcal D_{A_1}\times \mathcal D_{A\backslash A_1,B}^+.
\end{equation}
Clearly, $2|A_1|-1\geq 0$ and $2(|A|-|A_1|)+|B|-1\geq 0$.
Intuitively, this corresponds to the situation, where points in $D^\times$ labelled by $A_1$ collapse together to the origin.
\item[$ii)$] There is a subset $A_1$ of $A$, obeying $2\leq |A_1|\leq |A|$, such that
\begin{equation}\label{eq-diskb2}
\partial_{A_1}\mathcal D_{A,B}^+\cong \mathcal C_{A_1}\times \mathcal D_{A\backslash A_1\sqcup \{\bullet\},B}^+.
\end{equation}
We must impose $2|A_1|-3\geq 0$ and $2(|A|-|A_1|+1)+|B|-1\geq 0$.
Intuitively, this corresponds to the situation, where points in~$D^\times$ labelled by $A_1$ collapse together to a~single point in~$D^\times $ labelled by $\bullet$.
\item[$iii)$] Finally, there is a subset $A_1$ of $A$ and an ordered subset $B_1$ of successive elements of $B$,
obeying $0\leq |A_1|\leq |A|$ and $2\leq |B_1|\leq |B|$, such that
\begin{equation}\label{eq-diskb3}
\partial_{A_1,B_1}\mathcal D_{A,B}^+\cong \mathcal C_{A_1,B_1}^+\times
\mathcal D_{A\backslash A_1,B\backslash B_1\sqcup \{\bullet\}}^+.
\end{equation}
We impose $2|A_1|+|B_1|-2\geq 0$ and $2(|A|-|A_1|)+(|B|-|B_1|+1)-1\geq 0$.
Intuitively, this corresponds to the situation, where points in $D^\times$ labelled by $A_1$ and points in $S^1$ labelled by $B_1$ collapse together to a single point in $S^1$ labelled by $\bullet$.
\end{enumerate}

\subsection[An identification between compactified configuration spaces]{An identif\/ication between compactif\/ied conf\/iguration spaces}\label{ss-2-3}

We may use the action of $S^1$ to construct a section of $D_{A,B}^+$, namely we f\/ix one point $\circ$ in $S^1$ to~$1$.
This section is dif\/feomorphic, by means of the M\"obius transformation
\[
\psi\,:\,\mathcal H\sqcup\mathbb R\,\longrightarrow\,D\sqcup S^1\backslash\{\mathrm 1\};
\quad z\,\longmapsto\,\frac{z-\mathrm i}{z+\mathrm i}\,,
\]
where $D$ is the unit disk, to a smooth section of $C_{A\sqcup\{\bullet\},B\backslash\{\circ\}}^+$, given by
f\/ixing one point $\bullet$ in the complex upper half-plane $\mathcal H$ to $\mathrm i$ by means of the action of $G_2$.

Then, the compactif\/ied conf\/iguration space $\mathcal D_{A,B}^+$ can be identif\/ied with $\mathcal C_{A\sqcup\{\bullet\},B\backslash\{\circ\}}^+$,
and we observe that the cyclic order on the points in $S^1$ translates naturally into an order on the points
on the real axis $\mathbb R$.

We further consider the manifold $D_A$, and notice the identif\/ication $D_A\cong C_{A\sqcup\{\bullet\}}$: to be more
precise, by means of complex translation, we may put one point $\bullet$ in $C_{n+1}$ at the origin, and using
rescalings, one can put the remaining points in the punctured unit disk with boundary.
Analogously as before, the compactif\/ication $\mathcal D_A$ of $D_A$ can be identif\/ied with $\mathcal C_{A\sqcup\{\bullet\}}$.

We consequently identify the codimension $1$ boundary strata of $\mathcal D_{A,B}^+$ with those of $\mathcal C_{A\sqcup\{\bullet\},B\backslash\{\circ\}}^+\!\!$ (higher codimension can be worked out along the same lines very easily):
\begin{itemize}\itemsep=0pt
\item[$i)$] A boundary stratum as in \eqref{eq-diskb1} corresponds to the situation, where points labelled by
$A_1\sqcup\{\bullet\}$ collapse together to a single point in $\mathcal H$, which takes the r\^ole of the marked point~$\bul$.
\item[$ii)$] A boundary stratum as in \eqref{eq-diskb2} corresponds to the situation, where points labelled by $A_1$
collapse together to a single point in $\mathcal H$, which will not be the new marked point~$\bul$.
\item[$iii_1)$] A boundary stratum as in \eqref{eq-diskb3}, where $\circ\notin B_1$, corresponds to the situation, where points labelled
by $A_1\sqcup B_1$ collapse to a single point in $\mathbb{R}$, which will not be the new marked point~$\circ$.
\item[$iii_2)$] Finally, a boundary stratum as in \eqref{eq-diskb3}, where $\circ\in B_1$, corresponds to the situation, where points labelled
by the set $A\backslash A_1\sqcup\{\bullet\}\sqcup B\backslash B_1$ collapse to a single point in $\mathbb{R}$, which will be the new marked point~$\circ$.
\end{itemize}

\section[Explicit formulae for the formality morphism for cochains and chains]{Explicit formul\ae\ for the formality morphism\\ for cochains and chains}\label{s-3}

Here is a short review of the formul\ae\ we will need to construct the aforementioned $L_\infty$-quasi-isomorphisms $\mathcal U$ and $\mathcal S$.

\subsection[The $L_\infty$-quasi-isomorphism $\mathcal U$]{The $\boldsymbol{L_\infty}$-quasi-isomorphism $\boldsymbol{\mathcal U}$}\label{ss-3-1}

For the sake of simplicity, we denote by $[n]$, for a positive integer $n$, the set $\{1,\dots,n\}$.
For any pair of non-negative integers $(n,m)$, such that $2n+m-2\geq 0$, an {\em admissible graph} $\Gamma$ of
type $(n,m)$ is by def\/inition a directed graph with labels obeying the following requirements:
\begin{enumerate}\itemsep=0pt
\item[$i)$] The set of vertices $V_\Gamma$ is given by $[n]\sqcup [m]$; vertices labelled by integers in $[n]$,
resp.\ $[m]$, are called vertices of the f\/irst, resp.\ second type; further, the labelling of vertices of the
f\/irst type specif\/ies an order on them.
The set of vertices factorizes into $V_\Gamma=V_\Gamma^1\sqcup V_\Gamma^2$, where $V_\Gamma^1$, resp.\ $V_\Gamma^2$,
is the set of vertices of the f\/irst type, resp.\ second type.
\item[ii)] Every edge in $E_\Gamma$ starts at some vertex of the f\/irst type; there is at most one edge between
any two distinct vertices of $\Gamma$; no edge starts and ends at the same vertex.
\end{enumerate}
For a given vertex $v$ of $\Gamma$, we denote by $\mathrm{star}(v)$ the subset of $E_\Gamma$ of edges starting at
$v$: then, we assume that, for any vertex of the f\/irst type $v$ of $\Gamma$, the elements of $\mathrm{star}(v)$ are
labelled as $(e_v^1,\dots,e_v^{|\mathrm{star}(v)|})$.
By def\/inition, the {\em valence} of a vertex $v$ is the cardinality of the star of $v$.
The set of admissible graphs of type $(n,m)$ is denoted by $\mathcal G_{n,m}$.

We also need the following lemma, borrowed from~\cite{K}, to which we also refer for a more careful explanation
of the origins of the form $\omega$.
\begin{lemma}\label{l-angle}
There exists a smooth $1$-form $\omega$ on $\mathcal C_{2,0}$, with the following properties:
\begin{enumerate}\itemsep=0pt
\item[$i)$] The restriction of $\omega$ to the boundary stratum $\mathcal C_2=S^1$ equals the deRham differential of the angle function measured
in counterclockwise direction from the positive imaginary axis.
\item[$ii)$] The restriction of $\omega$ to $\mathcal C_{1,1}$, where the first point in the complex upper half-plane goes to the real axis, vanishes.
\end{enumerate}
\end{lemma}
For any pair of non-negative integers $(n,m)$, such that $2n+m-2\geq 0$, there are natural smooth projections from $\mathcal C_{n,m}^+$ onto
$\mathcal C_{2,0}$ (provided $n\geq 2$) or onto $\mathcal C_{1,1}$ (if $n,m\geq 1$), extending the natural projections on the open conf\/iguration spaces.

To an admissible graph $\Gamma$ of type $(n,m)$ is associated its Kontsevich's weight $W_\Gamma$ via
\begin{equation}\label{eq-wU}
W_\Gamma=\int_{\mathcal C_{n,m}^+} \bigwedge_{e\in E_\Gamma}\omega_e=\int_{\mathcal C_{n,m}^+} \omega_\Gamma,
\end{equation}
where, for an edge $e$ of $\Gamma$, $\omega_e$ denotes the pull-back of $\omega$ to $\mathcal C_{n,m}^+$ via the projection $\pi_e$ from
$\mathcal C_{n,m}^+$ onto $\mathcal C_{2,0}$, onto the pair of points labelled by the endpoints of $e$.
The labelling of $\Gamma$ specif\/ies an order of the forms $\omega_e$ in the above product.

To an admissible graph $\Gamma$ of type $(n,m)$, to $n$ polyvector f\/ields $\gamma_1,\dots,\gamma_n$ and to $m$ functions $f_1,\dots,f_m$ on $V$,
such that $|\mathrm{star}(k)|=|\gamma_k|+1$, $k=1,\dots,n$, we associate a function
\[
U_\Gamma(\gamma_1,\dots,\gamma_n)(f_1,\dots,f_m)
\]
by the following rule: to a vertex $v$ of the f\/irst type, resp.\ second type, we associate the polyvector f\/ield $\gamma_v$, resp.\ the function $f_v$.
For a function $I$ from $E_\Gamma$ to $[d]$, we associate to the vertex $v$ of the f\/irst type, resp.\ second type, the function
\[
\varphi_v^I=\gamma_v^{I(e_v^1),\dots,I(e_v^{|\mathrm{star}(v)|})},\qquad \text{resp.}\quad \varphi_v^I=f_v.
\]
with the same notations as before.
The function $I$ labels edges of $\Gamma$ by (standard) coordinates of $V$.
Then, we have the following assignment, for an admissible graph $\Gamma$ of type $(n,m)$:
\begin{equation}\label{eq-diffU}
U_\Gamma(\gamma_1,\dots,\gamma_n)(f_1,\dots,f_m)=\sum_{I:E_\Gamma\to [d]}\prod_{v\in V_\Gamma}\left(\prod_{e\in E_\Gamma:e=(*,v)}\partial_{I(e)}\right)\varphi_v^I.
\end{equation}
It is clear that $U_\Gamma$, with $\Gamma$ as above, maps $n$ polyvector f\/ields $\{\gamma_1,\dots,\gamma_n\}$ to a polydif\/ferential operator of (shifted)
degree $m-1$ by its very construction.

Finally, we def\/ine the $n$-th Taylor component $\mathcal U_n$ of Kontsevich's $L_\infty$-quasi-isomorphism $\mathcal U$ by combining (\ref{eq-wU})
and (\ref{eq-diffU}), namely
\begin{equation}\label{eq-taylorU}
\mathcal U_n=\sum_{m\geq 0}\sum_{\Gamma\in\mathcal G_{n,m}}W_\Gamma U_\Gamma.
\end{equation}
\begin{theorem}[Kontsevich]\label{t-kont}
The Taylor components \eqref{eq-taylorU} combine to an $L_\infty$-quasi-isomorphism
\[
\mathcal U:T_\mathrm{poly}^\bullet(V)\to \mathcal D_\mathrm{poly}^\bullet(V),
\]
of $L_\infty$-algebras, whose first order Taylor component reduces to the Hochschild--Kostant--Rosen\-berg quasi-isomorphism in cohomology.
\end{theorem}
For a proof of Theorem~\ref{t-kont}, we refer to~\cite{K}.

\subsection[The $L_\infty$-quasi-isomorphism $\mathcal S$]{The $\boldsymbol{L_\infty}$-quasi-isomorphism $\boldsymbol{\mathcal S}$}\label{ss-3-2}

An admissible graph of type $(n,m)$, for any two non-negative integers such that $2n+m-1\geq 0$, is a directed labelled graph $\Gamma$ as in
Subsection~\ref{ss-3-1}, with the only dif\/ference that there is a~special vertex, labelled by $0$, from which edges can only depart; the
vertex $0$ belongs neither to vertices of the f\/irst type nor of the second type.
The other requirements and notations from Subsection~\ref{ss-3-1} remain unaltered.
The set of admissible graphs of this kind of type $(n,m)$ is denoted by $\mathcal G_{n,m,0}$.

A (partial) counterpart of Lemma~\ref{l-angle} in this framework is the following lemma.

We def\/ine a smooth $1$-form on the conf\/iguration space $C_{3,0}$ via
\begin{equation*}
\varphi_D(p,q,r)=\varphi(q,r)-\varphi(q,p),
\end{equation*}
for any three pairwise distinct points $p$, $q$, $r$ in $\mathcal H\sqcup\mathbb R$: we then set $\omega_D=\mathrm d\varphi_D$.
\begin{lemma}\label{l-angleD}
The $1$-form $\omega_D$ extends to a smooth $1$-form on $\mathcal C_{3,0}$, with the
following properties:
\begin{itemize}\itemsep=0pt
\item[$i)$] its restriction to $\mathcal C_{2,1}$, when $q$ approaches the real axis, vanishes;
\item[$ii)$] its restriction to $\mathcal C_{2,0}\times\mathcal C_{1,1}$,
when $p$ and $q$ collapse together to the real axis,
equals $-\pi_1^*\omega$;
\item[$iii)$] its restriction to $\mathcal C_2\times\mathcal C_{2,0}$,
when $p$ and $q$ collapse together in the upper half-plane,
equals $\pi_2^*\omega-\pi_1^*\omega$;
\item[$iv)$] its restriction to $\mathcal C_{2,0}\times\mathcal C_{1,1}$ (resp.~$\mathcal C_2\times\mathcal C_{2,0}$),
when $p$ and $r$ collapse together to the real axis (resp.~in the upper half-plane),
vanishes;
\item[$v)$] its restriction to $\mathcal C_{2,0}\times\mathcal C_{1,1}$,
when $q$ and $r$ collapse together to the real axis,
equals $\pi_1^*\omega$;
\item[$vi)$] its restriction to $\mathcal C_2\times\mathcal C_{2,0}$,
when $q$ and $r$ collapse together in the upper half-plane,
equals $\pi_1^*\omega-\pi_2^*\omega$.
\end{itemize}
\end{lemma}

As above, we def\/ine a Shoikhet's weight associated to a graph without loop $\Gamma$ with
$m+n+1$ vertices labelled by $\mathcal{V}(\Gamma):=\{0,\dots,n,\overline{1},\dots,\overline{m}\}$.
To any edge $e=(i,j)\in\mathcal E(\Gamma)$, we associate a~smooth $1$-form $\omega_{D,e}$ on $\mathcal D_{n,m}^+$ by the following rules:
\begin{itemize}\itemsep=0pt
\item if neither $i$ nor $j$ lies in $\{0,\overline{1}\}$, then $\omega_{D,e}:=\pi_{(0,i,j)}^*\omega_D$, where
\begin{gather*}
\pi_{(0,i,j)}\,:\,\mathcal D_{n,m}^+\cong\mathcal C_{n+1,m-1}^+   \longrightarrow   \mathcal{C}_{3,0}, \quad
\big[(z_0,\dots,z_n,z_{\overline{2}},\dots,z_{\overline{m}})\big]    \longmapsto   \big[(z_0,z_i,z_j)\big] ;
\end{gather*}
\item if $i=0$ and $j\neq1$, then $\omega_{D,e}:=\pi_{(i,j)}^*\omega$, where
\begin{gather*}
\pi_{(i,j)}\,:\,\mathcal D_{n,m}^+\cong\mathcal C_{n+1,m-1}^+   \longrightarrow   \mathcal{C}_{2,0}, \quad
\big[(z_0,\dots,z_n,z_{\overline{2}},\dots,z_{\overline{m}})\big]  \longmapsto   \big[(z_i,z_j)\big] ;
\end{gather*}
\item if $j=1$ and $i\neq0$, then $\omega_{D,e}:={\rm p}_{(i,j)}^*\omega$, where
\begin{gather*}
{\rm p}_{(i,j)}\,:\,\mathcal D_{n,m}^+   \longrightarrow   \mathcal D_{1,1}\cong\mathcal{C}_{2,0}, \quad
\big[(z_1,\dots,z_n,z_{\overline{1}},\dots,z_{\overline{m}})\big]   \longmapsto   \big[(z_i,z_j)\big] ;
\end{gather*}
\item if $i=1$ or $j=0$ or $(i,j)=(0,1)$, then $\omega_{D,e}=0$.
\end{itemize}
Then, as above,
\begin{equation*}
\omega_{D,\Gamma}:=\bigwedge_{e\in \mathcal E(\Gamma)}\omega_{D,e}
\end{equation*}
def\/ines a dif\/ferential form on $\mathcal D_{n,m}^+$.
\begin{definition}
The Shoikhet weight $W_{D,\Gamma}$ of the directed graph $\Gamma$ is
\begin{equation*}
W_{D,\Gamma}:=\int_{\mathcal D_{n,m}^+} \omega_{D,\Gamma} .
\end{equation*}
\end{definition}

We consider an admissible graph in $\mathcal G_{n,m,0}$, such that $|\mathrm{star}(0)|=l$.
To $n$ polyvector f\/ields $\{\gamma_1,\dots,\gamma_n\}$ on $V$, such that $|\mathrm{star}(k)|=|\gamma_k|+1$, $k=1,\dots,n$, and to a Hochschild
chain $c=(a_0|a_1|\cdots|a_{m-1})$ of degree $-m+1$, we associate an $l$-form on $V$ (whose actual degree is $-l$, following the grading
introduced in~\cite{Sh}) via
\begin{equation*}
\langle\alpha,S_{\Gamma}(\gamma_1,\dots,\gamma_n,c)\rangle=
U_{\Gamma}(\alpha,\gamma_1,\dots,\gamma_n)(a_0,\dots,a_n)\,.
\end{equation*}

Finally, the $n$-th Taylor component $\mathcal S_n$ of the $L_\infty$-quasi-isomorphism $\mathcal S$ from the $L_\infty$-module $C_{-\bullet}(A,A)$
to the $L_\infty$-module $\Omega^{-\bullet}(V)$ over $T_\mathrm{poly}^\bullet(V)$ is given by
\begin{equation}\label{eq-taylorS}
\mathcal S_n=\sum_{m\geq 1}\sum_{\Gamma\in\mathcal G_{n,m,0}}W_{\Gamma,D}S_\Gamma.
\end{equation}
\begin{theorem}[Shoikhet]\label{t-shoikhet}
The Taylor components \eqref{eq-taylorS} combine to an $L_\infty$-quasi-isomorphism
\[
\mathcal S:C_{-\bullet}(A,A)\to \Omega^{-\bul}(V),
\]
of $L_\infty$-modules over $T_{\mathrm{poly}}^\bullet(V)$, whose $0$-th order Taylor component reduces to the
Hochschild--Kostant--Rosenberg quasi-isomorphism in homology.
\end{theorem}
We refer to~\cite{Sh} for the proof of Theorem~\ref{t-shoikhet}.

\section[Compatibility between the actions of polyvector fields on forms and Hochschild chains]{Compatibility between the actions of polyvector f\/ields\\ on forms and Hochschild chains}\label{s-4}

In this Section we sketch the proof of Theorem~\ref{t-capprod}, Section~\ref{s-1}, whose strategy owes to the homotopy argument used
in~\cite{K,MT}; for a more detailed version of the proof of Theorem~\ref{t-capprod} in an even more general case, we refer to~\cite{CR1}.

We must consider separately the case, where $\mathcal S_\gamma$ acts on Hochschild chains $i)$ of degree $m=0$ and $ii)$ of degree
$-m\leq -1$: as we will soon notice, the geometric aspects of the two cases are quite dif\/ferent.

\subsection[The space $\mathcal D_{1,1}$ and Hochschild chains of degree 0]{The space $\boldsymbol{\mathcal D_{1,1}}$ and Hochschild chains of degree 0}\label{ss-4-1}

By the def\/inition of the action $\cap$, if $c$ has degree $0$, the only Hochschild cochains acting on $c$ non-trivially must be functions,
in which case the action is simply multiplication on the right w.r.t.\ the product $\star$.

We consider the curve $\ell$ on the conf\/iguration space $\mathcal D_{1,1}$, with initial point on $\alpha$, and f\/inal point $b$, which
corresponds to the following embedding of the open unit interval into $D_{1,1}$:

\begin{figure}[h]
  \centerline{\includegraphics[scale=0.38]{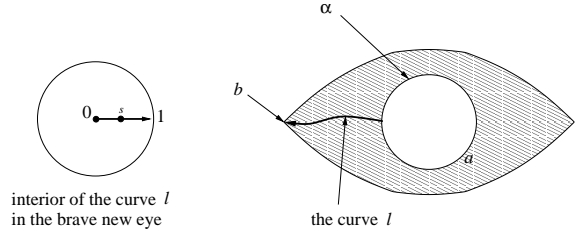}}
\caption{The curve $\ell$ in $\mathcal D_{1,1}$.}
\end{figure}

Since $\mathcal
D_{1,1}\cong\mathcal
C_2$, $\mathcal D_{1,1}$ coincides with
Kontsevich's eye: its ``pupil'' $\alpha$ represents the boundary stratum
$\mathcal D_2\times \mathcal D_{0,1}$ of codimension $1$, while the point $b$ represents the boundary stratum
$\mathcal C_{1,0}\times \mathcal C_{0,2}^+\times \mathcal D_{0,1}$ of codimension $2$, graphically
\begin{figure}[h]
  \centerline{\includegraphics[scale=0.35]{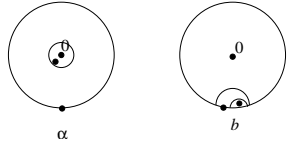}}
\caption{The boundary strata $\alpha$ and $b$ of $\mathcal D_{1,1}$.}
\end{figure}

The subset $\mathcal Y_{n,m}^+$ of $\mathcal D_{n,m}^+$, for $n\geq 1$, consisting of those conf\/igurations, whose projection onto~$\mathcal D_{1,1}$ belongs to the curve $\ell$, is a smooth submanifold with corners of $\mathcal D_{n,m}^+$ of codimension~$1$.
Pictorially, a typical conf\/iguration in $\mathcal Y_{n,1}$ looks like as follows:
\begin{figure}[h]
  \centerline{\includegraphics[scale=0.35]{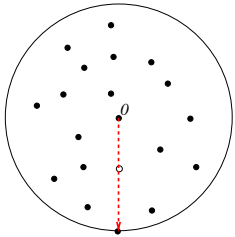}}
\caption{A typical conf\/iguration in $\mathcal Y_{n,1}$.}
\end{figure}

The dashed line represents the curve, along which the f\/irst point in $D^\times$ (labelled by $\circ$) moves, going from the origin to the unit circle.

We are interested in the boundary strata of $\mathcal Y_{n,1}$ of codimension $1$, which correspond to the chosen point on the pupil $\alpha$
and to the point $b$ of $\mathcal D_{1,1}$, namely
\begin{enumerate}\itemsep=0pt
\item[$i)$] conf\/igurations in $\mathcal D_{n,1}$, whose projection onto $\mathcal D_{1,1}$ is the initial point of the curve $\ell$
(the corresponding strata are denoted collectively by $\mathcal Y_{n,1}^0$);
\item[$ii)$] conf\/igurations in $\mathcal D_{n,1}$, whose projection onto $\mathcal D_{1,1}$ is the f\/inal point of the curve $\ell$
(the corresponding strata are denoted collectively by $\mathcal Y_{n,1}^1$);
\end{enumerate}

Let $\gamma$ be a solution of the Maurer--Cartan equation as in (\ref{eq-MCel}), $\alpha$ a polyvector f\/ield on $V$ of degree $-1$ and
$c=a_0$ a Hochschild chain of degree $0$ for the algebra $A$.
\begin{proposition}\label{p-cap}
For $\gamma$, $\alpha$ and $c$ as above, the following identities hold true:
\begin{gather}
\label{eq-cap-0}\alpha\cap \mathcal S_\gamma(c)=\sum_{n\geq 0}\frac{1}{n!}\sum_{\Gamma\in\mathcal G_{n+1,1,0}}W_{D,\Gamma}^0 \mathcal S_\Gamma(\alpha,\underset{n}{\underbrace{\gamma,\dots,\gamma}},c),\\
\label{eq-cap-1}\mathcal S_\gamma(\mathcal U_\gamma(\alpha)\cap c)=\sum_{n\geq 0}\frac{1}{n!}\sum_{\Gamma\in\mathcal G_{n+1,1,0}}W_{D,\Gamma}^1 S_\Gamma(\alpha,\underset{n}{\underbrace{\gamma,\dots,\gamma}},c),
\end{gather}
where the weights $W_{D,\Gamma}^i$, $i=0,1$, are defined via
\[
W_\Gamma^i=\int_{\mathcal Y_{n,1}^i}\omega_{D,\Gamma}.
\]
\end{proposition}
\begin{proof}
The proof of (\ref{eq-cap-0}) and (\ref{eq-cap-1}) relies mainly on the evaluation of the weights $W_{D,\Gamma}^i$, $i=0,1$: we only give a sketch of such evaluations, referring to~\cite{CR1} for all details.

By construction, all admissible graphs $\Gamma$ appearing in the expressions on the right-hand side of (\ref{eq-cap-0}) and (\ref{eq-cap-1})
have all vertices of the f\/irst type of valence $2$, except the vertex labelled by $1$, which in this case has valence $0$.

First, for any admissible graph $\Gamma$ in $\mathcal G_{n+1,1,0}$, the weight $W_{D,\Gamma}^0$ vanishes, if $1$ has at least one incoming edge: namely, if $1$ has one arrow coming from $0$, then the integral vanishes, since the angle form is the derivative of a (locally) constant function.
Further, Kontsevich's Vanishing Lemma~\cite[Lemma 6.4]{K} applies to the remaining cases, whence only the case matters, where $0$ and $1$ collapse together, and then Lemma~\ref{l-angleD} does the job.
Otherwise, the identity
\begin{equation}\label{eq-wS-0}
W_{D,\Gamma}^0=W_{D,\Gamma_0}
\end{equation}
holds true, where $\Gamma_0$ is the admissible graph in $\mathcal G_{n,1}$, obtained from $\Gamma$ by collapsing the vertices~$0$ and $1$.
Here is a graphical representation of a general component $Z$ of $\mathcal Y_{n+1,1}^0$
\begin{figure}[h]
  \centerline{\includegraphics[scale=0.35]{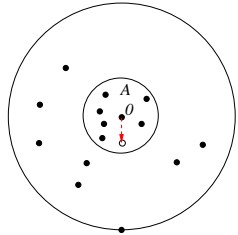}}
\caption{A typical conf\/iguration in $\mathcal Y_{n,1}^0$.}
\end{figure}

Second, for an admissible graph $\Gamma$ in $\mathcal G_{n+1,1,0}$, the weight $W_{D,\Gamma}^1$, restricted to a component $Z$ of $\mathcal Y_{n+1,1}^1$ of the form
\[
Z\cong \mathcal C_{A_1,0}\times\mathcal C_{A_2,2}^+\times\mathcal D_{A_3,1},\qquad 1\leq |A_1|,
\]
vanishes, unless there are no edges outgoing from $A_1$ or from $A_2$, in which cases the weight factorizes as
\begin{equation}\label{eq-wS-1}
W_{D,\Gamma}^1|_Z=W_{\Gamma_1}W_{\Gamma_2}W_{D,\Gamma_3},
\end{equation}
and $\Gamma_1$ is the admissible subgraph of $\Gamma$, whose vertices
of the first type are labelled by $A_1$, $\Gamma_2$~is the graph, whose
vertices are labelled by $A_2\sqcup \{1,2\}$, and obtained by
collapsing~$\Gamma_1$ to the vertex~$2$, and $\Gamma_3$ is the graph,
whose vertices are labelled by $A_3\sqcup \{1\}$, obtained by
collapsing~$\Gamma_2$ to the vertex~$1$.

\newpage

Graphically, a typical conf\/iguration in the component $Z$ looks like as
\begin{figure}[h!]
  \centerline{\includegraphics[scale=0.35]{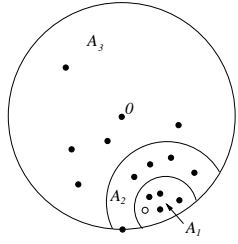}}
\caption{A typical conf\/iguration in $\mathcal Y_{n,1}^1$.}
\end{figure}

First of all, if there is an edge e.g.\ from $A_1$ to $A_2$, the corresponding contribution vanishes by means of Lemma~\ref{l-angleD}, $i)$; the same argument implies the claim in all other cases.
We observe that this also implies that $\Gamma_i$, $i=1,2,3$, is admissible.

Further, the f\/irst two factors in the factorization of the weight $W_{D,\Gamma}^1$ reduce to usual Kontsevich's weights, again in virtue of Lemma~\ref{l-angleD}, $v)$.
\end{proof}

The curve $\ell$ in $\mathcal D_{1,1}$ ``interpolates'' between $\mathcal Y_{n+1,1}^0$ and $\mathcal Y_{n+1,1}^1$: we may evaluate weights of admissible graphs in $\mathcal G_{n+1,1,0}$ on the remaining boundary strata of codimension $1$ of $\mathcal Y_{n+1,1}$, and, by means of Stokes' Theorem, this implies that (\ref{eq-cap-0}) and (\ref{eq-cap-1}) coincide at the level of cohomology; for a complete discussion of the corresponding homotopy formula, we refer to~\cite{CR1}.

\subsection[The space $\mathcal D_{1,2}^+$ and Hochschild chains of higher degree]{The space $\boldsymbol{\mathcal D_{1,2}^+}$ and Hochschild chains of higher degree}\label{ss-4-2}

We prove now Theorem~\ref{t-capprod} in the case, where (\ref{eq-tangh}) is applied to Hochschild chains $c$ of higher (negative) degree.

The open unit square in $\mathbb C$ can be embedded in the open conf\/iguration space $D_{1,2}^+$ via $(s,t)\mapsto \left[\left(s,1,e^{2\pi \mathrm i t}\right)\right]$, where the square brackets denote equivalence classes w.r.t.\ the action of $S^1$; we may take a possible closure $\sigma$ of it in the compactif\/ication $\mathcal D_{1,2}^+$ as follows:
\begin{figure}[h]
  \centerline{\includegraphics[scale=0.35]{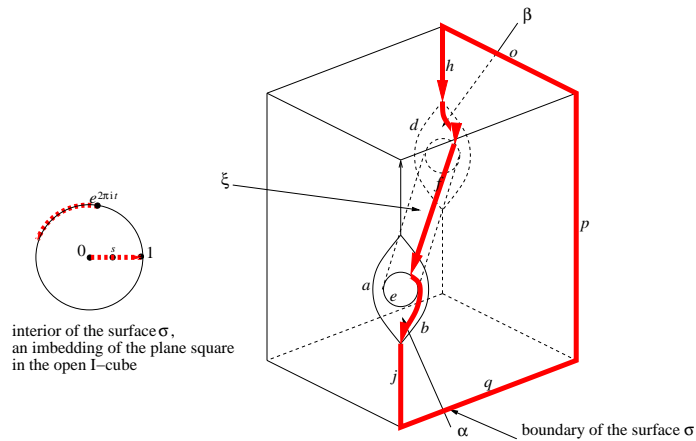}}
\caption{The boundary of $\sigma$ in $\mathcal D_{1,2}^+$.}
\end{figure}

Topologically, $\mathcal D_{1,2}^+$ is a ``cube with two eyes'', or I-cube: we will be mostly interested, in the forthcoming discussion, in the boundary stratum $\xi$ of codimension $1$, which is $\mathcal D_1\times \mathcal D_{0,2}^+$, and in the boundary strata $o$, $q$ of codimension $2$, which are $\mathcal C_{1,0}\times\mathcal C_{0,2}^+\times \mathcal D_{0,2}^+$ and $\mathcal C_{1,1}\times\mathcal C_{0,2}^+\times\mathcal D_{0,1}$ respectively; graphically
\begin{figure}[h]
  \centerline{\includegraphics[scale=0.35]{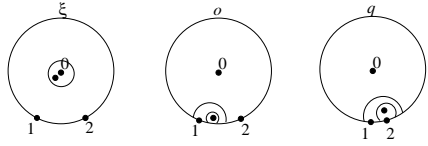}}
\caption{The boundary strata $\xi$, $o$ and $q$ of $\mathcal D_{1,2}^+$.}
\end{figure}

We observe that the straight line on the boundary stratum $\xi$ corresponds to the edge $\{s=0\}$, the boundary stratum $o$ corresponds to the edge $\{s=1\}$ and the boundary stratum $q$ corresponds to (a way of imbedding) the point $(1,0)$.

For any two positive integers $(m,n)$, with $n\geq 1$ and $m\geq 2$, the subset $\mathcal Y_{n,m}^+$ of those conf\/i\-gu\-rations in $\mathcal D_{n,m}^+$, whose projection onto $\mathcal D_{1,2}^+$ belongs to $\sigma$, is a smooth, orientable submanifold with corners of $\mathcal D_{n,m}^+$ of codimension $1$.
Graphically,
\begin{figure}[h]
  \centerline{\includegraphics[scale=0.35]{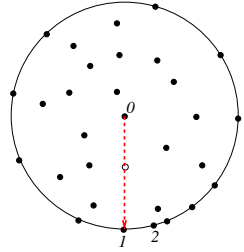}}
\caption{A typical conf\/iguration in $\mathcal Y_{n,m}^+$.}
\end{figure}

We will need the boundary strata of $\mathcal Y_{n,m}^+$ of codimension $1$, corresponding to conf\/igurations in~$\mathcal D_{n,m}^+$, whose projection onto the I-cube is in the boundary of $\sigma$ (collectively denoted by~$\mathcal Y_{n,m}^{+,\partial\sigma}$).
More precisely, $\mathcal Y_{n,m}^{+,\partial\sigma}$ factorizes into eight dif\/ferent types, denoted by $\mathcal Y_{n,m}^{+,x}$, and we will consider only $x$ to be the straight line on $\xi$ or $o$ and~$q$.

We consider a solution $\gamma$ of the Maurer--Cartan equation as in (\ref{eq-MCel}), a polyvector f\/ield $\alpha$ on~$V$, and a Hochschild chain $c=(a_0|\cdots|a_m)$ in $A$, $m\geq 1$.
\begin{proposition}\label{p-cap-I}
For $\gamma$, $\alpha$ and $c$ of degree $-m\leq -1$ as above, the following identities hold true:
\begin{equation}\label{eq-cap-xi}\alpha\cap \mathcal S_\gamma(c)=\sum_{n\geq 0}\frac{1}{n!}\sum_{\Gamma\in\mathcal G_{n+1,m+1,0}}W_{D,\Gamma}^\xi \mathcal S_\Gamma(\alpha,\underset{n}{\underbrace{\gamma,\dots,\gamma}},c),
\end{equation}
and
\begin{equation}\label{eq-cap-oq}
\mathcal S_\gamma(\mathcal U_\gamma(\alpha)\cap c) =\begin{cases}
\sum_{n\geq 0}\frac{1}{n!}\sum_{\Gamma\in\mathcal G_{n+1,m+1,0}}W_{D,\Gamma}^o S_\Gamma(\alpha,\underset{n}{\underbrace{\gamma,\dots,\gamma}},c),& |\alpha|=-1,\\
\sum_{n\geq 0}\frac{1}{n!}\sum_{\Gamma\in\mathcal G_{n+1,m+1,0}}W_{D,\Gamma}^q S_\Gamma(\alpha,\underset{n}{\underbrace{\gamma,\dots,\gamma}},c),& |\alpha|\geq 0,
\end{cases}
\end{equation}
where the weights $W_{D,\Gamma}^x$, $x=\xi,o,q$, are defined via
\[
W_{D,\Gamma}^x=\int_{\mathcal Y_{n+1,m}^{+,x}}\omega_{D,\Gamma}.
\]
\end{proposition}
\begin{proof}
The proof follows along the same lines of the proof of Proposition~\ref{p-cap}, with some due changes; once again, we refer to~\cite{CR1} for the complete proofs, while here we make some necessary comments on the degrees, which hold true in this particular situation.

We observe that, for any admissible graph $\Gamma$ in $\mathcal G_{n+1,m,0}$, the weight $W_{D,\Gamma}^\xi$ vanishes, if the vertex $1$ has at least one incoming edge, by the very same arguments sketched in the proof of Proposition~\ref{p-cap}.
Otherwise, the identity
\[
W_{D,\Gamma}^\xi=W_{D,\Gamma_0}
\]
holds true, where $\Gamma_0$ in $\mathcal G_{n,m,0}$ is obtained from $\Gamma$ by collapsing the vertices $0$ and $1$: this generalizes (\ref{eq-wS-0}) in the proof of Proposition~\ref{p-cap}, and the proof uses almost the same arguments.

On the other hand, a general component $Z$ of the boundary strata of $\mathcal Y_{n,m}^{+,o}$, resp.\ $\mathcal Y_{n,m}^{+,q}$ has the explicit form
\begin{gather}
\label{eq-b-o}Z \cong \mathcal C_{A_1,0}^+\times\mathcal C_{A_2,B_2}^+\times\mathcal D_{A_3,B_3}^+,\quad \text{resp.}\\
\label{eq-b-q}Z \cong \mathcal C_{A_1,B_1}^+\times\mathcal C_{A_2,B_2}^+\times\mathcal D_{A_3,B_3}^+,
\end{gather}
where $A_i$, $i=1,2,3$, are disjoint subsets of $[n]$, with $1\leq |A_1|\leq n$, $0\leq |A_2|\leq n$, $0\leq |A_3|\leq n$, with $n=|A_1|+|A_2|+|A_3|$, and $B_i$, $i=1,2,3$, are disjoint ordered subsets of $[m]$ of successive elements, such that $1\leq |B_1|\leq m$, $2\leq |B_2|\leq m$, $1\leq |B_3|\leq m$, and $m=|B_1|+|B_2|+|B_3|$.
Pictorially,
\begin{figure}[h]
  \centerline{\includegraphics[scale=0.35]{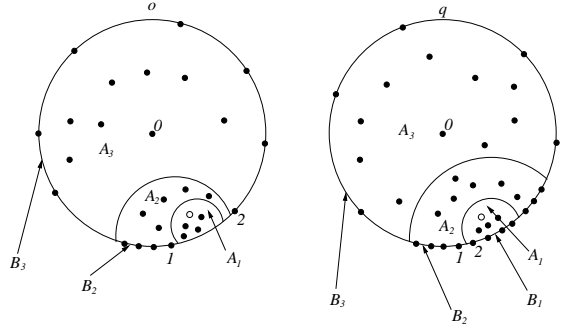}}
\caption{Typical components of $\mathcal Y_{n,m}^{+,o}$ and $\mathcal Y_{n,m}^{+,q}$.}
\end{figure}

We consider a component $Z$ of $\mathcal Y_{n+1,m}^o$ as in (\ref{eq-b-o}), resp.\ of $\mathcal Y_{n+1,m}^q$ as in (\ref{eq-b-q}), and for an admissible graph $\Gamma$ as before, we denote by $\Gamma_1$ the admissible subgraph of
$\Gamma$, whose vertices are labelled by $A_1\sqcup B_1$, by $\Gamma_2$
the graph, whose vertices are labelled by $A_2\sqcup B_2$, and obtained
by collapsing $\Gamma_1$ to a single vertex, and by $\Gamma_3$ the graph,
whose vertices are labelled by $A_3\sqcup B_3$, obtained by collapsing
$\Gamma_2$ to a single vertex.

The weight $W_{D,\Gamma}^o$, restricted to a component $Z$ of $\mathcal Y_{n+1,m}^{+,o}$ as above vanishes, unless there are no edges going from $A_1$ to $A_2$ or $A_3$, or from $A_2$ to $A_3$, $|B_2|=2$, and $\alpha$ has degree $-1$.
In fact, the weight $W_{D,\Gamma}^o$ factorizes as
\[
W_{D,\Gamma}^o|_Z=W_{\Gamma_1}W_{\Gamma_2}W_{D,\Gamma_3},
\]
with $B_1=\varnothing$.
In particular, the third weight on the right-hand side is non-trivial only if $2|A_1|-2$ equals the degree of the integrand, which, by the above reasonings, is precisely $2(|A_1|-1)+|\alpha|+1$ (since all vertices of the f\/irst type in $A_1$ are $2$-valent except the f\/irst vertex): therefore, the integral is non-trivial only if $|\alpha|=-1$.

The weight $W_{D,\Gamma}^q$, restricted to a component $Z$ of $\mathcal Y_{n+1,m}^{+,q}$ as above vanishes, unless there are no edges going from $A_1$ to $A_2$ or $A_3$, or from $A_2$ to $A_3$, $|B_2|=2$, and $|B_1|=|\alpha|+1$ (which implies that $\alpha$ has degree bigger or equal than $0$).
In such cases, the weight $W_{D,\Gamma}^o$ factorizes as
\[
W_{D,\Gamma}^o|_Z=W_{\Gamma_1}W_{\Gamma_2}W_{D,\Gamma_3},
\]
with $1\leq |B_1|$.
Dimensional reasons as in the case of a component $Z$ of $\mathcal Y_{n+1,m}^{+,q}$ force the degree of $\alpha$ to be bigger or equal than $0$: namely, the degree of the third integrand is $2(|A_1|-1)+|\alpha|+1$, and the integral is non-trivial only if it equals $2|A_1|+|B_1|-2$.
Since $1\leq |B_1|$, the non-triviality condition forces $|\alpha|=|B_1|-1$, whence the claim.

This result obviously generalizes (\ref{eq-wS-1}) in the proof of Proposition~\ref{p-cap}, and its proof is the same as the proof of (\ref{eq-wS-1}).
\end{proof}
The surface $\sigma$ ``interpolates'' between the boundary strata $\mathcal Y_{n+1,m}^{+,\xi}$, $\mathcal Y_{n+1,m}^{+,o}$ and $\mathcal Y_{n+1,m}^{+,q}$: the ``interpolation'' in this situation is of course more complicated than the one in Subsection~\ref{ss-4-1}, since we have to keep track of the boundary of $\sigma$.
In fact, the weighted sums as in~(\ref{eq-cap-xi}) and~(\ref{eq-cap-oq}), where we integrate over boundary strata of $\mathcal Y_{n+1,m}^{+,x}$, $x=h,j,p$, vanish; the weighted sums over the remaining boundary strata of $\mathcal Y_{n+1,m}^+$ can be also explicitly evaluated, and, by means of Stokes' Theorem, they produce a homotopy formula, proving that the left-hand sides of (\ref{eq-cap-xi}) and (\ref{eq-cap-oq}) coincide at the level of cohomology.

\section[Duflo isomorphism on (co)invariants revisited]{Duf\/lo isomorphism on (co)invariants revisited}\label{s-5}

Let $\mg$ be a f\/inite dimensional Lie algebra over $\mathbb C$.

First of all, the {\em (modified) Duflo element}~\cite{Du} is def\/ined via
\[
J:=\det\Big(\frac{e^{\mathrm{ad}/2}-e^{-\mathrm{ad}/2}}{\mathrm{ad}}\Big)\in\widehat{\S}(\mg^*)^{\mg} .
\]

We have a morphism of $\mg$-modules
\[
\mathcal D\,:=\,{\rm sym}\circ(J^{1/2}\cdot):\S(\mg)\,\longrightarrow\,\U(\mg),
\]
where $\U(\mg)$ is the Universal Enveloping Algebra of $\mg$, and $\rm{sym}$ denotes the usual symmetrization map from the symmetric algebra $\S(\mg)$ to $\U(\mg)$.

The following result generalizes to coinvariants the well-known Duf\/lo isomorphism~\cite{Du}.
\begin{theorem}\label{t-chdi}
The map $\mathcal D$ restricts to an isomorphism of algebras
$\S(\mg)^\mg\,\tilde\longrightarrow\,\U(\mg)^\mg=\mathcal Z\big(\U(\mg)\big)$ on invariants,
and induces an isomorphism of $\S(\mg)^\mg$-modules
$\S(\mg)_\mg\,\tilde\longrightarrow\,\U(\mg)_\mg=\mathcal{A}\big(\U(\mg)\big)$ on coinvariants.
\end{theorem}
Here $\mathcal{Z}(B)$ denotes the center of an algebra $B$, and $\mathcal A(B)=B/[B,B]$ its abelianization.

We sketch here a proof of Theorem~\ref{t-chdi} in the spirit of the approach of \cite{K,PT} to the original Duf\/lo isomorphism.

We consider the Kirillov--Kostant--Souriau Poisson bivector $\pi$ on $\mg^*$ and associated product~$\star$.
Since the product $\star$ obeys $x\star y-y\star x=[x,y]_{\mg}$\footnote{We may set $\hbar=1$, as the Poisson
structure is linear.}, viewing $x,y\in\mg$ as linear functions on $\mg^*$, there is an algebra morphism
\[
\mathcal{I}\,:\,\U(\mg)\longrightarrow\big(\S(\mg),\star\big)\,;\, x\longmapsto x .
\]
The map $\mathcal{I}^{-1}\circ\mathcal{U}_{\gamma}$ induces an algebra isomorphism by means of Theorem~\ref{t-cupprod}, Section~\ref{s-1},
\[
\S(\mg)^{\mg}=\mathcal{Z}\big(\S(\mg),\{,\}\big)\longrightarrow\mathcal{Z}\big(\U(\mg)\big)=\U(\mg)^\mg
\]
while, dually, $\mathcal{S}_{\gamma}\circ\mathcal{I}$ induces
an isomorphism of $\S(\mg)^\mg$-modules by means of Theorem~\ref{t-capprod}, Section~\ref{s-1},
\[
\U(\mg)_\mg=\mathcal{A}\big(\U(\mg)\big)\longrightarrow\mathcal{A}\big(\S(\mg),\{,\}\big)=\S(\mg)_{\mg}.
\]

In~\cite{Sh2}, the restriction of $\mathcal{U}_{\gamma}$ to $\S(\mg)$ has been shown to be the identity, which, coupled with Kontsevich's discussion \cite[Section 8]{K}, implies that
$\mathcal{I}^{-1}\circ\mathcal{U}_{\gamma}=\mathcal{I}^{-1}=\mathcal{D}:\S(\mg)\to\U(\mg)$.

Shoikhet's proof of the fact that the restriction of $\mathcal U_{\gamma}$ to functions is the identity
can be summarized as follows: for any function $f\in\S(\mg)$, $\mathcal{U}_{\gamma}(f)$ is expressed only
{\it via} so-called inner wheels, whose weights vanish by the main result of~\cite{Sh2}.
Pictorially, an inner wheel looks like as follows:

\begin{figure}[h!!]
  \centerline{\includegraphics[scale=0.39]{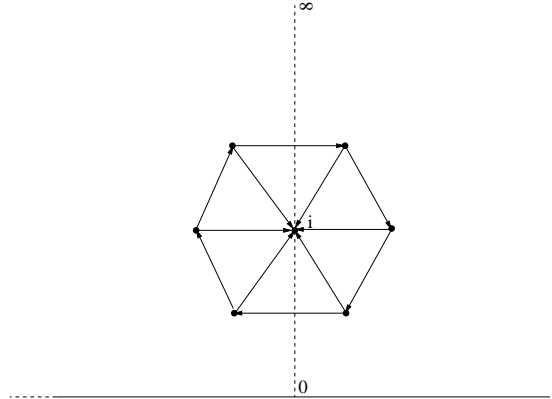}}
\caption{An inner wheel.}
\end{figure}

Theorem~\ref{t-chdi} follows then from the following
\begin{proposition}\label{p-vanwheel}
The restriction of $\mathcal{S}_{\gamma}$ on $\S(\mg)$ coincides with the identity map.
\end{proposition}
\begin{proof}
By the arguments of~\cite[Paragraph 3.6.1]{Sh}, the only admissible graphs contributing non-trivially to $\mathcal{S}_{\gamma}$ are of the form
\begin{figure}[h]
  \centerline{\includegraphics[scale=0.30]{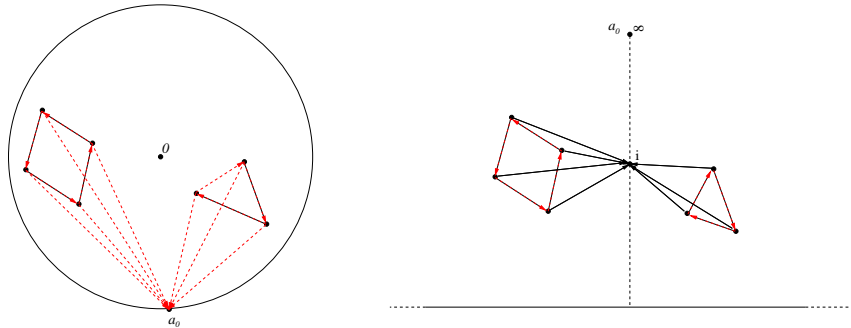}}
\caption{A typical admissible graph in $\mathcal S_{\gamma}$ on functions and its counterpart in the upper half-plane.} \label{fig11}
\end{figure}

Further, we use a correspondence between weights of admissible graphs in $\mathcal G_{n,1,0}$ and
admissible graphs in $\mathcal G_{n+1,0}$, by means of the M\"obius transformation
\[
\psi\,:\,\mathcal{H}\sqcup\mathbb{R}\,\longrightarrow\, D\sqcup S^1\backslash\{1\};\quad z\,\longmapsto\,\frac{z-i}{z+i} ,
\]
which induces (see Subsection \ref{ss-2-3}) isomorphisms $\mathcal{C}_{n+1,0}\,\tilde\longrightarrow\,\mathcal{D}_{n,1}$, to prove that
$\mathcal S_\gamma={\rm id}$.

First of all, given on $\mathcal D_{1,1}$ the $1$-form $\omega_D$ as in Lemma \ref{l-angleD}, then the $1$-form
$\omega:=\psi^*\omega_D$ on $\mathcal C_{2,0}$ is a dif\/ference of usual Kontsevich's angle forms as in Lemma \ref{l-angle}. Then the weight $W_{D,\Gamma}$ of
an admissible graph $\Gamma$ in $\mathcal G_{n,1,0}$ is pulled-back to a weight $W_{\Gamma'}$, $\Gamma'$ being
admissible in $\mathcal G_{n+1,0}$:
the vertex $0$ is mapped to a vertex of the f\/irst type, while the only vertex of the second type is mapped to $\infty$
in the complex upper half-plane.
The factors of $\omega_{D,\Gamma}$ are pulled-back to $i)$ usual forms $\omega_e$,
whenever $e$ is an edge from some vertex (of the f\/irst and of the second type) to $1$, and the ``new'' $e$ is now an
edge from the (inverse image of the) starting point to $\mathrm i$, and $ii)$ dif\/ferences between $\omega_e$ and
$\omega_{e(\mathrm i)}$, if $e$ is an edge between two vertices (of the f\/irst and second type, neither of which is $0$)
and the new edge $e$ connects the (inverse images of the) endpoints, while $e(\mathrm i)$ is an edge, whose starting
point is the starting point of $e$ and whose endpoint is $\mathrm i$.
We have used here the arguments exposed in Appendix 1 of~\cite{CR1}, to which we refer for a more detailed discussion, as well as for the properties of $\omega_D$.
We also used the fact that the form $\omega$ vanishes when its f\/inal point goes to inf\/inity, which f\/inally implies the above correspondence.

We observe that dashed arrows denote forms $\omega_{D,e}$ on $\mathcal D_{n,1}$ in the left-most graph of Fig.~\ref{fig11}, while we have
used plain, resp.\ dashed, arrows to denote forms $\omega_e$ on $\mathcal C_{n+1,0}$, resp.\ dif\/ferences of such forms,
in the right-most graph.

Finally, we use the following graphical calculus for replacing dashed edges by plain ones:
\begin{figure}[h]
  \centerline{\includegraphics[scale=0.5]{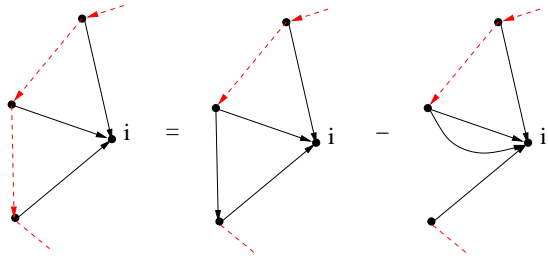}}
\caption{Replacing dashed edges by plain ones in an inner wheel.}
\end{figure}

The second graph on the right-hand side has a double edge, whence its weight vanishes.
Thus, the Shoikhet weight of the left-most graph in Fig.~\ref{fig11} equals the usual Kontsevich weight of the
right-most one (i.e.~with all edges turned into black ones), which is zero by~\cite{Sh2}.
\end{proof}

We f\/inally observe that Proposition~\ref{p-vanwheel} has been proved in Subsubsection 3.6.2 of~\cite{Sh} under the assumption of the validity of Conjecture 3.5.3.1, which is implied by Theorem~\ref{t-capprod}: our proof, on the contrary, is purely based on the main result of~\cite{Sh2} and on the properties of the forms~$\omega_D$.

\section{Conclusion}
In this paper, we have proved Shoikhet's conjecture using conf\/iguration space integrals: it is worthwhile noticing that we did not exploit all boundary strata of the I-cube, which replaces in our homotopy argument Kontsevich's eye (namely, we made use only of the boundary stratum~$\xi$ of codimension $1$ and of the boundary strata $o$ and $q$ of codimension~$2$).
In fact, Shoikhet's conjecture can be viewed as a special case of a more general result, which involves a Maurer--Cartan element of a more general shape, i.e.\ a sum of polyvector f\/ields of dif\/ferent degrees, to which, via Kontsevich's formality, corresponds also a sum of polydif\/ferential operators, also of dif\/ferent degrees: accordingly, all boundary strata of the surface $\sigma$ in the I-cube contribute to the proof of this more general result~\cite{CR1}.

\subsubsection*{Acknowledgement}

We thank Giovanni Felder for useful discussions and comments.
The research of D.C.~(on leave of absence from Universit\'e Lyon 1) is fully supported by the European
Union thanks to a Marie Curie Intra-European Fellowship (contract number MEIF-CT-2007-042212).

\pdfbookmark[1]{References}{ref}
\LastPageEnding

\end{document}